\documentclass[acmsmall,screen]{MapleTrans}
\usepackage[utf8]{inputenc}

\usepackage{tikz}
\usetikzlibrary{decorations.markings,patterns}
\usetikzlibrary{arrows.meta, decorations.markings}
\usepackage{pgfplots}
\usepackage[normalem]{ulem}

\usepackage{enumitem}
\usepackage{etoolbox}
\usepackage[utf8]{inputenc}
\usepackage{breqn}
\usepackage{maple}
\usepackage{xfrac}
\usepackage{xcolor}
\definecolor{mygreen}{rgb}{0,0.6,0}
\definecolor{mygray}{rgb}{0.5,0.5,0.5}
\definecolor{mymauve}{rgb}{0.58,0,0.82}
\definecolor{altblue}{rgb}{0.0,0.6,1.0}
\definecolor{lstbg}{cmyk}{0.05, 0.01, 0, 0}
\definecolor{morebluish}{cmyk}{0.06,0.04,0,0}
\usepackage{listings}
\input{listings-maple-definition.sty}
\setcopyright{acmlicensed}
\copyrightyear{2025}
\acmYear{2025}
\acmDOI{10.5206/mt.vViI.XXXXX}

\acmJournal{MAPLETRANS}
\acmVolume{X}
\acmNumber{Y}
\acmArticle{Z}
\acmMonth{0}

\begin{document}

\title[MT Template 1.0]{A Pedagogical Introduction to the Unified Transform Method: The Heat Equation on a Finite Interval}

\author{Athanasios Paraskevopoulos}
\authornotemark[1]
\email{at.paraskevopoulos@aua.gr}
\affiliation{%
  \institution{Mathematics Research Center, Academy of Athens/ Department of Natural Resources Development and Agricultural Engineering, Agricultural University of
Athens }
  \city{Athens}
  \country{Greece}
  \postcode{}
}

\renewcommand{\shortauthors}{Paraskevopoulos}

\begin{abstract}
This paper presents a detailed application of the Unified Transform Method (Fokas method) to the one-dimensional heat equation on $[0,1]$ with Dirichlet boundary conditions. The analysis formulates the Initial-Boundary Value Problem and derives an integral representation of the solution via a generalised spatial Fourier transform with complex spectral parameter $\lambda \in \mathbb{C}$, yielding the Global Relation --- an algebraic identity coupling the initial datum, prescribed boundary values, and unknown Neumann data. The unknowns are eliminated by exploiting the symmetry $\lambda \mapsto -\lambda$, reducing the solution to a contour integral over $\partial D^+$. An explicit evaluation is carried out for exponential initial datum $u_0(x)=e^{-x}$ and Dirichlet conditions $g_0(t)=\cos(t)$, $h_0(t)=e^{-1}\cos(t)$. The integral representation is analysed in the complex plane, with emphasis on exponential decay and analyticity, providing rigorous justification for contour deformation via Cauchy's Theorem and Jordan's Lemma. Numerical implementation in Maple uses a trapezoidal contour parametrisation ensuring exponential decay along each segment; the solution over $x\in[0,1]$, $t\in[0,2\pi]$ matches prescribed data to machine precision. The results confirm the analytical and numerical efficacy of the Unified Transform for classical parabolic problems and illustrate how rigorous contour analysis yields stable, accurate solutions.
\end{abstract}


\maketitle
\section{Introduction}
The analysis of Initial-Boundary Value Problems (IBVPs) for linear evolutionary partial differential equations has a long and well-developed history, rooted in the classical methods of separation of variables and the Fourier series. While these techniques are effective for canonical domains and standard boundary conditions, they lack the flexibility to handle more general settings in a unified manner. A significant departure from this tradition came in 1997 with the introduction of the Unified Transform Method, due to Fokas \cite{fokas1997unified}, subsequently developed and extended in \cite{flyer2008numerical, fokas2008unified}. This method furnishes a systematic framework for a broad class of IBVPs by synthesising spectral analysis with the theory of functions of a complex variable, yielding integral representations of the solution that are both analytically tractable and numerically implementable.

The distinguishing feature of the Fokas method, relative to classical approaches, lies in its treatment of the boundary data. In the standard Fourier approach, the unknown solution is expanded in a basis adapted to the domain, and the boundary conditions are imposed term by term. The Fokas method instead proceeds by applying a transform pair — a generalisation of the classical Fourier transform — directly to the governing equation, without a priori restricting the spectral variable \(\lambda\)
 to the real line. This produces a relation, valid for all \(\lambda \in \mathbb{C}\), which couples the initial data, the prescribed boundary values, and the unknown boundary data in a single equation. This object, referred to throughout as the Global Relation, plays a central role in the method: by exploiting the symmetries of the problem in the complex \(\lambda\)-plane, it is possible to solve for the unknown boundary data algebraically, without recourse to the governing equation itself.
 
The inversion of the transform and solution representation are carried out by deforming integration contours in the complex plane. This step requires careful justification, relying on classical results from complex analysis — principally Cauchy's Integral Theorem and Jordan's Lemma \cite{ablowitz2003complex} — which ensure that contributions from arcs at infinity vanish and that contours may be deformed freely within regions of analyticity. A key feature of the method is that the final representation involves integrals over contours along which the integrand is exponentially decaying, making the formula well-suited to numerical evaluation.

Despite the analytical depth and broad applicability of the Fokas method, it remains largely absent from standard graduate curricula in applied mathematics, and its exposition in the literature is often directed at specialist audiences. The present paper is intended to address this gap by providing a rigorous yet pedagogically accessible account of the method, suitable for students with a working knowledge of Fourier analysis and introductory complex variable theory. The problem treated throughout is the heat equation on the finite interval [0,1] subject to Dirichlet boundary conditions — a setting simple enough to admit a complete and explicit treatment, yet sufficiently structured to illustrate every essential step of the method, from the derivation of the Global Relation to the deformation of contours and the final numerical implementation in \textsc{Maple}.

The paper is organised as follows. In \autoref{2}, we formulate the IBVP precisely, apply the Fokas transform, derive the Global Relation, and obtain an integral representation of the solution expressed entirely in terms of prescribed data. In \autoref{3}, we specialise to an explicit example with exponential initial data and cosine boundary conditions. We simplify the integral representation analytically, show that one of the integral contributions vanishes by Cauchy's theorem, and numerically evaluate the remaining contour integral in \textsc{Maple} using a trapezoidal contour parametrisation. The section concludes with a three-dimensional visualisation of the solution over \(x \in [0,1]\) and \(t \in [0,2\pi]\).

\section{Problem formulation}\label{2}
We study the one-dimensional heat conduction equation on the finite interval, subject to Dirichlet  boundary conditions \cite{chatziafratis2025fokas}:
\begin{subequations}\label{bvp}
\begin{alignat}{3}
u_t &= u_{xx}, && 0 < x < 1,\; t > 0  \\
u(x,0) &= u_0(x), && 0 < x < 1  \\
u(0,t) &= g_0(t), u(1,t)=h_0(t) && t > 0 \label{bvp4}
\end{alignat}
\end{subequations}

\begin{figure}[H]
\centering
\begin{tikzpicture}[scale=1.2]
  \fill[gray!30] (0,0) rectangle (4,2);
  \draw[->, thick] (-0.3,0) -- (4.5,0) node[anchor=west] {\( x \)};
  \draw[->, thick] (0,-0.3) -- (0,3.5) node[anchor=south] {\( t \)};
  \draw[blue, thick] (0,2) -- (4,2);
  \draw[blue, thick] (4,0) -- (4,2);
  
\end{tikzpicture}
\caption{Domain of interest for the heat equation}
\end{figure}
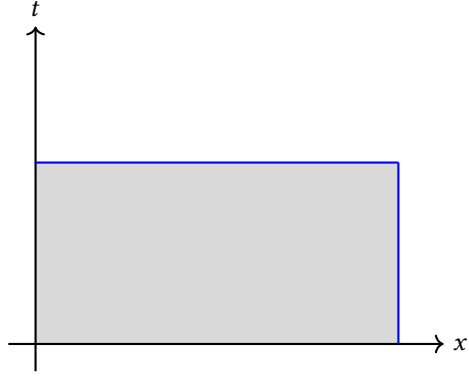

We apply the Fourier transform in space to reduce the PDE to an ODE in the spectral variable \(\lambda\). Then we solve the ODE and take the inverse transform to recover \( u(x,t) \).

 We assume that there exists a function \(u(x,t)\), \(0<x<1\) and \(t>0\) that satisfies the PDE and the initial conditions.

Applying the Fourier transform of both sides of the PDE \((1\text{a})\):
\[\mathcal{F}\{u_{t}\}=\mathcal{F}\{u_{xx}\}\]

to the heat equation and integrating by parts twice yields the following first-order ODE in \(t\)

    \begin{align*}
        \frac{\partial \hat{u}}{\partial t}&+\lambda^2\hat{u}(\lambda,t)=(u_x(1,t)+i\lambda u(1,t))e^{-i\lambda}-(u_x(0,t)+i\lambda u(0,t))
  \end{align*}
  
This is a first-order inhomogeneous ODE in time for \(\hat{u}(\lambda,t)\). We use the integrating factor method.

The integrating factor is:
\(\mu(t)=e^{\int p(t) \ dt}=e^{\int \lambda^2 \ dt}=e^{\lambda^2t}\).

We multiply all terms of the differential equation by the integrating factor.

\begin{align*}
   (e^{\lambda^2t}\hat{u}(\lambda,t))'&=e^{\lambda^2t-\lambda i}(u_x(1,t)+i\lambda u(1,t)) -(u_x(0,t)+i\lambda  u(0,t))e^{\lambda^2t}
\end{align*}

Integrating both sides with respect to \(t\) gives

\begin{align}
e^{\lambda^2t}\hat{u}(\lambda,t)-\hat{u}_0(\lambda)=e^{-\lambda i} \int_{0}^{t}u_x(1,t)e^{\lambda^2\tau} d\tau+e^{-\lambda i} \lambda i \int_{0}^{t}u(1,t)e^{\lambda^2\tau} d\tau-\int_{0}^{t} u_x(0,\tau) e^{\lambda^2\tau} \:d\tau - i\lambda \int_{0}^{t} u(0,\tau) e^{\lambda^2\tau} \:d\tau  \label{sec_formulation}
\end{align}

Now we define the integral transform as 

\begin{equation*}
\tilde{g}_{j}(\lambda,t) =  \int_0^t g(s) e^{\lambda s} \, ds,  \, t>0, \, j=0,1,\, \lambda \in \mathbb{C}
\end{equation*}

with \(g_1(t)=u_x(0,t)\), \(g_0(t)=u(0,t),\, t>0\)
and

\begin{align*}
\tilde{h}_{j}(\lambda,t) =  \int_0^t h(s) e^{\lambda s} \, ds,  \, t>0, \, j=0,1,\, \lambda \in \mathbb{C}
\end{align*}

with \(h_1(t)=u_x(1,t)\), \(h_0(t)=u(1,t),\, t>0\).

In this notation, the integrated equation takes the form

\begin{align}
    \label{3}
   e^{\lambda^2t}\hat{u}(\lambda,t)-\hat{u}_0(\lambda)=-\tilde{{g}}_1(\lambda^2,t)-i \lambda \tilde{g}_0(\lambda^2,t)+ e^{-\lambda i} \Big(\tilde{h}_1(\lambda^2,t)+i\lambda  \tilde{h}_0(\lambda^2,t)\Big) 
\end{align}

This is called the Global Relation (G.R). This equation is valid for all \(\lambda \in \mathbb{C}\) and \(t>0\).

 Inverting the Fourier transform formally gives:

\begin{equation} \label{4}
\begin{aligned}
    u(x,t) = \frac{1}{2\pi} \int_{-\infty}^{\infty} e^{i\lambda x - \lambda^2 t} \hat{u}_0(\lambda)\,d\lambda 
    &- \frac{1}{2\pi}\int_{-\infty}^{\infty} e^{i\lambda x - \lambda^2 t} 
    \Big[\tilde{g}_1(\lambda^2,t) + i\lambda\,\tilde{g}_0(\lambda^2,t)\Big]\,d\lambda \\
    &+ \frac{1}{2\pi}\int_{-\infty}^{\infty} e^{-i\lambda(1-x) - \lambda^2 t} 
    \Big[\tilde{h}_1(\lambda^2,t) + i\lambda\,\tilde{h}_0(\lambda^2,t)\Big]\,d\lambda
\end{aligned}
\end{equation}

	   For the first two integrals, that is the formula for the heat equation
on the half-line, where the contour \(\partial D^{+}\) is the boundary of the domain \(D^{+}\) defined by

\[D^{+}=\{\operatorname{Im}(\lambda) \geq 0, \ \operatorname{Re}(\lambda^2)<0\}\]

For the third integral, the only difference is that now, instead of the factor \(e^{i\lambda x}\), there exists the factor \(e^{i \lambda(x-1)}\).

Let the \(\lambda=\lambda_R+i\lambda_I \in \mathbb{C}\) then 
:
\begin{itemize}
    \item \(|e^{i \lambda_R(x-1)}|=1\)
    \item The growth or decay comes from \(e^{\lambda_I (1-x)}\)
    \begin{itemize}
        \item if \(\lambda_I<0\), it decays as \(x \rightarrow 1\).
        \item if \(\lambda_I>0\), this term grows as \(x \rightarrow 1\).
    \end{itemize}
\end{itemize}

We define the region \( D^- \), which is the mirror image of \( D^+ \) with respect to the real axis. Specifically:
\[
D^- = \left\{ \lambda \in \mathbb{C} : \operatorname{Im}(\lambda) \leq 0,\ \operatorname{Re}(\lambda^2) < 0 \right\}.
\]

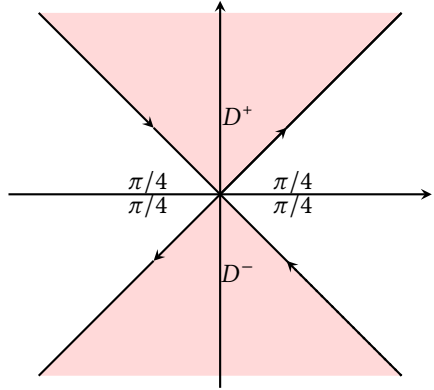
\begin{figure}[H]
\centering
\begin{tikzpicture}[scale=.8, >=stealth]

\fill[red!15] (-3,3) -- (0,0) -- (3,3) -- cycle;

\fill[red!15] (-3,-3) -- (0,0) -- (3,-3) -- cycle;

\draw[->, thick] (-3.5, 0) -- (3.5, 0) node[below right] {};
\draw[->, thick] (0, -3.2) -- (0, 3.2) node[above] {};

\draw[thick, ->] (-3,3) -- (-1.1,1.1);
\draw[thick] (-1.1,1.1) -- (0,0);
\draw[thick, ->] (0,0) -- (1.1,1.1);
\draw[thick] (1.1,1.1) -- (3,3);

\draw[ thick, ->] (3,3) -- (-1.1,-1.1);
\draw[ thick] (-1.1,-1.1) -- (-3,-3);
\draw[dotted, thick, ->] (3,-3) -- (1.1,-1.1);
\draw[ thick] (3,-3) -- (0,0);

\node at (-1.2,0.2) {\(\pi/4\)};
\node at (1.2,0.2) {\(\pi/4\)};
\node at (-1.2,-0.2) {\(\pi/4\)};
\node at (1.2,-0.2) {\(\pi/4\)};

\node at (0.3,1.3) {\(D^+\)};
\node at (0.3,-1.3) {\(D^-\)};

\end{tikzpicture}
\caption{Domains \(D^+=\{\operatorname{Im}(\lambda) \geq 0, \ \operatorname{Re}(\lambda^2)<0\}\) and \(D^{-} = \{\operatorname{Im}(\lambda) \leq 0, \ \operatorname{Re}(\lambda^2)<0\}\), in the upper and lower half-planes, respectively. }
\end{figure}

\subsection{First integral:}
 \[\int_{-\infty}^{\infty} e^{i \lambda x
    -\lambda^2t} \hat{u}_0(\lambda)\:d\lambda\]

 This is the standard Fourier Inversion of the initial condition. If \(\hat{u}(\lambda)\) is the Fourier transform of a function defined on \(x>0\) and is sufficiently regular (e.g., rapidly decaying), this integral converges absolutely on the real line — no contour deformation is needed. This term is handled in the same manner as in standard PDE Fourier methods. No unknown boundary data is involved.

\subsection{Second integral:}

Suppose that only \(g_1(t)=u_x(0,t)\) is unknown, while \(g_0(t)=u(0,t)\) is prescribed. The second integral in the solution representation is

\begin{align*}
\int_{-\infty}^{\infty} e^{i \lambda x
    -\lambda^2t} (\tilde{{g}}_1(\lambda^2,t)+ i \lambda \tilde{g}_0(\lambda^2,t))\:d\lambda
\end{align*}

The result of the integral along this path \(C^+\) is
\begin{align*}
   &\oint_{C^+} e^{i \lambda x
    -\lambda^2t} (\tilde{{g}}_1(\lambda^2,t)+ i \lambda \tilde{g}_0(\lambda^2,t))\:d\lambda= \int_{-R}^{R} e^{i \lambda x
    -\lambda^2t} (\tilde{{g}}_1(\lambda^2,t)+ i \lambda \tilde{g}_0(\lambda^2,t))\:d\lambda\\&+ \int_{\gamma_1} e^{i \lambda x
    -\lambda^2t} (\tilde{{g}}_1(\lambda^2,t)+ i \lambda \tilde{g}_0(\lambda^2,t))\:d\lambda+\int_{\gamma_2} e^{i \lambda x
    -\lambda^2t} (\tilde{{g}}_1(\lambda^2,t)+ i \lambda \tilde{g}_0(\lambda^2,t))\:d\lambda \\
    &+\int_{\gamma_3} e^{i \lambda x
    -\lambda^2t} (\tilde{{g}}_1(\lambda^2,t)+ i \lambda \tilde{g}_0(\lambda^2,t))\:d\lambda + \int_{\gamma_4} e^{i \lambda x
    -\lambda^2t} (\tilde{{g}}_1(\lambda^2,t)+ i \lambda \tilde{g}_0(\lambda^2,t))\:d\lambda
\end{align*}

We want to examine the case where  \(R \rightarrow \infty\) since these are the bounds of our integral. Firstly, note that because it has no singularities within the contour \(C^+\),
\[\oint_{C^+} e^{i \lambda x
    -\lambda^2t} (\tilde{{g}}_1(\lambda^2,t)+ i \lambda \tilde{g}_0(\lambda^2,t))\:d\lambda=0\]

Now, we evaluate each integral in the contour decomposition. Consider the integral over \(\gamma_1\)
\begin{align*}
\int_{\gamma_1} e^{i \lambda x
    -\lambda^2t} (\tilde{{g}}_1(\lambda^2,t)+ i \lambda \tilde{g}_0(\lambda^2,t))\:d\lambda=\int_{\gamma_1}e^{i\lambda x} f(\lambda)\: d \lambda
    \end{align*}

where we define \begin{align*}
     f(\lambda)&=e^{-\lambda^2t} (\tilde{{g}}_1(\lambda^2,t)+ i \lambda \tilde{g}_0(\lambda^2,t))\\
    &=e^{-\lambda^2t}\int_{0}^{t}e^{\lambda^2s}g_1(s)\:ds+e^{-\lambda^2t}i \lambda \int_{0}^{t}e^{\lambda^2s}g_0(s)\:ds
\end{align*}

Under the polar representation of complex numbers,  \(\lambda=R e^{i \theta}\), \(0\leq \theta\leq \frac{\pi}{4}\) and \(d\lambda=Rie^{i\theta} \ d \theta\), so

\[
\int_{\gamma_1}e^{i\lambda x} f(\lambda)\: d \lambda=\int_{0}^{\frac{\pi}{4}} e^{ixRe^{i\theta}} f(Re^{i\theta}) Rie^{i\theta} \ d\theta
\]

An application of Jordan's Lemma \cite{ablowitz2003complex}

    while the function \(f\) is bounded along the arc \(\gamma_1\)    and as
\(R \rightarrow \infty , e^{-xR\sin{\theta}} \rightarrow 0\) and \(e^{-R^2\cos{2\theta}} \rightarrow 0\). 

So we have that 
\begin{equation}
\int_{\gamma_1} e^{i \lambda x
    -\lambda^2t} (\tilde{{g}}_1(\lambda^2,t)+ i \lambda \tilde{g}_0(\lambda^2,t))\:d\lambda=0
\end{equation}

Similarly, we can prove that 
\begin{equation}
   \int_{\gamma_2} e^{i \lambda x
    -\lambda^2t} (\tilde{{g}}_1(\lambda^2,t)+ i \lambda \tilde{g}_0(\lambda^2,t))\:d\lambda=0
\end{equation}

Thus,
 \begin{align*}
      &-\int_{-\infty}^{\infty} e^{i \lambda x
    -\lambda^2t} (\tilde{{g}}_1(\lambda^2,t)+ i \lambda \tilde{g}_0(\lambda^2,t))\:d\lambda\\&=\int_{\gamma_3} e^{i \lambda x
    -\lambda^2t} (\tilde{{g}}_1(\lambda^2,t)+ i \lambda \tilde{g}_0(\lambda^2,t))\:d\lambda + \int_{\gamma_4} e^{i \lambda x
    -\lambda^2t} (\tilde{{g}}_1(\lambda^2,t)+ i \lambda \tilde{g}_0(\lambda^2,t))\:d\lambda
\end{align*}

\subsection{Third integral:}
Similarly \(h_1(t)=u_x(1,t)\) is unknown, while \(h_0(t)=u(1,t)\) is prescribed. The third integral in the solution representation is

\[\int_{-\infty}^{\infty} e^{i \lambda (x-1)-\lambda^2t} \Big(\tilde{h}_1(\lambda^2,t)+i\lambda  \tilde{h}_0(\lambda^2,t)\Big)\]

The result of the integral along this path \(C^{-}\) is
\begin{align*}
    &\oint_{C^-} e^{i \lambda (x-1)-\lambda^2t} (\tilde{{h}}_1(\lambda^2,t)+ i \lambda \tilde{h}_0(\lambda^2,t))\:d\lambda= -\int_{-R}^{R}e^{i \lambda (x-1)-\lambda^2t} (\tilde{{h}}_1(\lambda^2,t)+ i \lambda \tilde{h}_0(\lambda^2,t))\:d\lambda\\&- \int_{\gamma_5} e^{i \lambda (x-1)-\lambda^2t} (\tilde{{h}}_1(\lambda^2,t)- i \lambda \tilde{h}_0(\lambda^2,t))\:d\lambda -\int_{\gamma_6} e^{i \lambda (x-1)-\lambda^2t} (\tilde{{h}}_1(\lambda^2,t)+ i \lambda \tilde{h}_0(\lambda^2,t))\:d\lambda \\
    &-\int_{\gamma_7} e^{i \lambda (x-1)-\lambda^2t} (\tilde{{h}}_1(\lambda^2,t)+ i \lambda \tilde{h}_0(\lambda^2,t))\:d\lambda -\int_{\gamma_8} e^{i \lambda (x-1)-\lambda^2t} (\tilde{{h}}_1(\lambda^2,t)+ i \lambda \tilde{h}_0(\lambda^2,t))\:d\lambda
\end{align*}

By applying Cauchy's Theorem and Jordan's Lemma, we obtain:

\begin{align*}
   \int_{-R}^{R}e^{i \lambda (x-1)-\lambda^2t} (\tilde{{h}}_1(\lambda^2,t)+ i \lambda \tilde{h}_0(\lambda^2,t))\:d\lambda=-\int_{\gamma_6} e^{i \lambda (x-1)-\lambda^2t} (\tilde{{h}}_1(\lambda^2,t)+ i \lambda \tilde{h}_0(\lambda^2,t))\:d\lambda \\
    -\int_{\gamma_7} e^{i \lambda (x-1)-\lambda^2t} (\tilde{{h}}_1(\lambda^2,t)+ i \lambda \tilde{h}_0(\lambda^2,t))\:d\lambda 
\end{align*}

So the formula \(\eqref{4}\) can be written as
\begin{equation}\label{7}
    \begin{aligned}
    u(x,t)&=\frac{1}{2\pi} \int_{-\infty}^{\infty} e^{i \lambda x
    -\lambda^2t} \hat{u}_0(\lambda)\:d\lambda-\frac{1}{2\pi}\int_{\partial D^{+}}e^{i \lambda x
    -\lambda^2t} (\tilde{{g}}_1(\lambda^2,t)+ i \lambda \tilde{g}_0(\lambda^2,t))\:d\lambda\\ &-\frac{1}{2\pi}\int_{\partial D^{+}}e^{i \lambda (1-x)-\lambda^2t}(-\tilde{{h}}_1(\lambda^2,t)+ i \lambda \tilde{h}_0(\lambda^2,t))\:d\lambda  
\end{aligned}
\end{equation}

Now, we know from G.R. \(\eqref{3}\), that the functions \(\hat{u}_{0}, \tilde{g}_0,\tilde{h}_0\) are known but the functions \(\hat{u} , \tilde{g}_1\) are unknown.

We define the function \(G\) as
\begin{equation}
G(\lambda,t)=\tilde{u}_0(\lambda)-i \lambda\tilde{g}_0(\lambda^2,t)+e^{-\lambda i}i \lambda \tilde{h}_0(\lambda^2,t)
\end{equation}

so the G.R \(\eqref{3}\) can be writen in the form of

\begin{equation}\label{9}
    e^{\lambda^2t} \hat{u}(\lambda,t)=G(\lambda,t)-\tilde{g}_1(\lambda^2,t)+e^{-i\lambda}\tilde{h}_1(\lambda^2,t) \end{equation}

Applying the transformation \(\lambda \mapsto - \lambda\) to \(\eqref{9}\) we obtain:

\begin{equation}\label{10}
    e^{\lambda^2t} \hat{u}(-\lambda,t)=G(-\lambda,t)-\tilde{g}_1(\lambda^2,t)+e^{i\lambda}\tilde{h}_1(\lambda^2,t) \end{equation}

Notice: \(\tilde{g}_1,\tilde{h}_1\) are even in \(\lambda\) (they depend on \(\lambda^2\)), so they stay the same under
\(\lambda \mapsto - \lambda\)

Now we have a system of two equations \(\eqref{9},\eqref{10}\) in two unknowns \(\tilde{g}_1,\tilde{h}_1\)

We multiply the equation \(\eqref{9}\) by \( e^{i\lambda } \):

\begin{equation}\label{11}
e^{i\lambda } e^{\lambda^2 t} \hat{u}(\lambda, t)
= e^{i\lambda } G(\lambda, t) - e^{i\lambda } \tilde{g}_1(\lambda^2,t) + \tilde{h}_1(\lambda^2,t) 
\end{equation}

Then we multiply equation \(\eqref{10}\) by \( -e^{-i\lambda } \):

\begin{equation}\label{12}
- e^{-i\lambda } e^{\lambda^2 t} \hat{u}(-\lambda, t)
= - e^{-i\lambda } G(-\lambda, t) + e^{-i\lambda } \tilde{g}_1(\lambda^2,t) - \tilde{h}_1(\lambda^2,t) 
\end{equation}

So we add the equations \(\eqref{11}\) and \(\eqref{12}\):

\begin{align*}
&e^{i\lambda } e^{\lambda^2 t} \hat{u}(\lambda, t)
- e^{-i\lambda } e^{\lambda^2 t} \hat{u}(-\lambda, t) \\
&= e^{i\lambda } G(\lambda, t)
- e^{-i\lambda } G(-\lambda, t)
+ \left( -e^{i\lambda } + e^{-i\lambda } \right) \tilde{g}_1
\end{align*}

Solving for \( \tilde{g}_1 \):

\begin{equation}
\begin{aligned}
\tilde{g}_1 &=
\frac{
e^{i\lambda } G(\lambda, t)
- e^{-i\lambda } G(-\lambda, t)
+ e^{\lambda^2 t} \left[
 e^{-i\lambda } \hat{u}(-\lambda, t)- e^{
 i\lambda } \hat{u}(\lambda, t)\right]
}{
e^{i\lambda } - e^{-i\lambda }
}\\
&= \frac{
e^{i\lambda } G(\lambda, t)
- e^{-i\lambda } G(-\lambda, t)
}{
\Delta \lambda}+ e^{\lambda^2 t} \hat{g}_1 
\end{aligned}
\end{equation}

 where \(\Delta \lambda=e^{i \lambda}-e^{-i\lambda}\) and \(\hat{g}_1=-\frac{1}{\Delta \lambda}\left[
 e^{i\lambda } \hat{u}(\lambda, t)- e^{
 -i\lambda } \hat{u}(-\lambda, t)\right]\)

Similarly,we add \(\eqref{9},\eqref{10}\) and  solving for \(\tilde{h}_1\):
\begin{equation}
\begin{aligned}
\tilde{h}_1(\lambda^2, t) &=
\frac{\left[ G(\lambda, t) - G(-\lambda, t) \right]+
e^{\lambda^2 t} \left[ \hat{u}(-\lambda, t) - \hat{u}(\lambda, t) \right]
}{
 e^{i\lambda }-e^{-i\lambda }
}\\&=\frac{\left[ G(\lambda, t) - G(-\lambda, t) \right]
}{
 \Delta \lambda
}+ e^{\lambda^2t} \hat{h}_1
\end{aligned}
\end{equation}

where \(\hat{h}_1=-\frac{1}{\Delta \lambda}[ \hat{u}(\lambda, t) - \hat{u}(-\lambda, t) ]\).

We next substitute \( \tilde{g}_1 \) and \( \tilde{h}_1 \) in \(\eqref{7}\).

We can observe that the terms involving \(\hat{u}(\pm \lambda,t)\) are unknown explicitly. To complete the solution, we must show that these terms do not contribute to the integral representation. 

The term in \(\tilde{g}_1\) involves the following contribution from  \(\hat{u}(\pm \lambda,t)\)
\begin{equation} 
 \frac{
 e^{-i\lambda } \hat{u}(-\lambda, t)
- e^{i\lambda } \hat{u}(\lambda, t)
}{
e^{i\lambda } - e^{-i\lambda }
}  \end{equation}

Now, recall the definition of the Fourier transform
\begin{align*}
\hat{u}(\lambda,t) = \int_0^L e^{-i\lambda x} u(x,t)\, dx, \,\,\hat{u}(-\lambda,t) = \int_0^L e^{i\lambda x} u(x,t)\, dx
\end{align*}

The relevant contribution to the solution integral is \begin{align*}
    \int_{\partial D^{+}}e^{i\lambda x} \frac{e^{-i\lambda}}{e^{i\lambda}-e^{-i\lambda}}\hat{u}(-\lambda,t)\:d\lambda
\end{align*}

We examine the behavior of the integrand
\[\frac{e^{-i\lambda}}{e^{i\lambda}-e^{-i\lambda}}\sim \frac{e^{-i\lambda}}{-e^{-i\lambda}}\sim -1,\quad |\lambda|\rightarrow\infty\]

We now show that \(\hat{u}(-\lambda,t)\) as \(|\lambda|\rightarrow\infty\)
\begin{align*}
    \hat{u}(-\lambda,t)&=\Big[\frac{e^{i\lambda x}u(x,t)}{i \lambda}\Big]_{x=0}^{x=1}-\frac{1}{i\lambda}\int_{0}^{1} e^{i \lambda x} \frac{\partial u}{\partial x}(x,t)\:dx\\
    &=\frac{u(1,t)e^{i\lambda}-u(0,t)}{i \lambda}-O(\frac{1}{\lambda})
\end{align*}

which is clearly bounded as \(\lambda \rightarrow \infty\) with \(Im(\lambda)\geq 0\).

Also \[\int_{\partial D^{+}} e^{i\lambda x} \frac{e^{i\lambda}}{e^{i\lambda}-e^{-i\lambda}}\hat{u}(\lambda,t)\:d\lambda\sim \int_{\partial D^{+}} e^{i\lambda x} e^{2i\lambda}\hat{u}(\lambda,t)\:d\lambda\]
if \(\lambda_I>0\) the term \(e^{2\lambda i}\) decays 
\begin{align*}
    \hat{u}(\lambda,t)&=\Big[\frac{e^{-i\lambda x}u(x,t)}{-i \lambda}\Big]_{x=0}^{x=1}+\frac{1}{i\lambda}\int_{0}^{1} e^{-i \lambda x} \frac{\partial u}{\partial x}(x,t)\:dx\\
    &=\frac{u(0,t)-u(1,t)e^{-i\lambda}}{i \lambda}+O(\frac{1}{\lambda})
\end{align*}

which is clearly bounded as \(\lambda \rightarrow \infty\).

Similarly, the term in \(\tilde{h}_1\).
We have 
\begin{align*}
    \int_{\partial D^{+}} e^{i\lambda(1-x)} \frac{\hat{u}(-\lambda,t)}{e^{i\lambda}-e^{-i\lambda}}\:d\lambda=   \int_{\partial D^{+}} e^{-i\lambda x}\frac{e^{i\lambda}}{e^{i\lambda}-e^{-i\lambda}} \hat{u}(-\lambda,t)\:d\lambda
\end{align*}
Thus, when integrating along \(\partial D^{+}\):
\begin{itemize}
    \item \(e^{-i\lambda x}\) decays for \(Im(\lambda)\geq 0\)
    \item  the factor \(\frac{e^{i\lambda}}{e^{i\lambda}-e^{-i\lambda}}\sim 1 \) as \(|\lambda|\rightarrow \infty\)
\end{itemize}

So we obtain

\begin{equation}
\oint_{\partial D^{+}} e^{i\lambda (1-x)}
\frac{
 \hat{u}(-\lambda, t)  -\hat{u}(\lambda, t)  
}{
e^{i\lambda }-e^{-i\lambda } 
}=0 \end{equation}

Thus, the formula \(\eqref{7}\) can be written as 
\begin{align*}
    u(x,t)&=\frac{1}{2\pi} \int_{-\infty}^{\infty} e^{i \lambda x
    -\lambda^2t} \hat{u}_0(\lambda)\:d\lambda\\
    &-\frac{1}{2\pi}\int_{\partial D^{+}} e^{ -\lambda^2t}  \Bigg\{e^{i \lambda x
   } \Big[\frac{
e^{i\lambda } G(\lambda, t)
- e^{-i\lambda } G(-\lambda, t)
}{
\Delta \lambda} + i \lambda \tilde{g}_0(\lambda^2,t)\Big]\\ 
&+e^{i \lambda (1-x)}\Big[-\frac{\left[ G(\lambda, t) - G(-\lambda, t) \right]
}{
 \Delta \lambda
}+ i \lambda \tilde{h}_0(\lambda^2,t)\Big]\Bigg\} \:d\lambda 
\end{align*}

We recall that \(G(\lambda,t)=\hat{u}_0(\lambda)-i \lambda\tilde{g}_0(\lambda^2,t)+e^{-\lambda i}i \lambda \tilde{h}_0(\lambda^2,t)\). We obtain

\begin{align*}
    u(x,t)&=\frac{1}{2\pi} \int_{-\infty}^{\infty} e^{i \lambda x
    -\lambda^2t} \hat{u}_0(\lambda)\:d\lambda\\
    &-\frac{1}{2\pi}\int_{\partial D^{+}} \frac{e^{ -\lambda^2t}}{\Delta \lambda}  \Bigg\{e^{i \lambda x
   } \Big[
e^{i\lambda } G(\lambda, t)
- e^{-i\lambda } G(-\lambda, t)
 + i \lambda (e^{i \lambda}-e^{-i\lambda})\tilde{g}_0(\lambda^2,t)\Big]\\ 
&+e^{i \lambda (1-x)}\Big[-\left[ G(\lambda, t) - G(-\lambda, t) \right]
+ i \lambda (e^{i \lambda}-e^{-i\lambda})\tilde{h}_0(\lambda^2,t)\Big]\Bigg\} \:d\lambda 
\end{align*}

The first big term 
\begin{align*}
  e^{i \lambda x
   } \Big[
e^{i\lambda } G(\lambda, t)
- e^{-i\lambda } G(-\lambda, t)
 + i \lambda (e^{i \lambda}-e^{-i\lambda})\tilde{g}_0(\lambda^2,t)\Big]&= e^{i \lambda x
   } \Big[
e^{i\lambda }(\hat{u}_0(\lambda)-i \lambda\tilde{g}_0(\lambda^2,t)+e^{-\lambda i}i \lambda \tilde{h}_0(\lambda^2,t))\\
&- e^{-i\lambda } (\hat{u}_0(-\lambda)+i \lambda\tilde{g}_0(\lambda^2,t)-e^{\lambda i}i \lambda \tilde{h}_0(\lambda^2,t))\\
& + i \lambda (e^{i \lambda}-e^{-i\lambda})\tilde{g}_0(\lambda^2,t)\Big] \\
&=e^{i\lambda x}[e^{i\lambda}\hat{u}_0(\lambda)-e^{-i\lambda} \hat{u}_0(-\lambda)-i\lambda(e^{i\lambda}+e^{-i\lambda})\tilde{g}_0+2i\lambda\tilde{h}_0\\
&+ i \lambda (e^{i \lambda}-e^{-i\lambda})\tilde{g}_0(\lambda^2,t)]\\
&=e^{i\lambda(x+1)}\hat{u}_0(\lambda)-e^{i\lambda(x-1)} \hat{u}_0(-\lambda)-2i\lambda e^{i \lambda(x-1)}\tilde{g}_0(\lambda^2,t)\\
&+2i\lambda e^{i\lambda x}\tilde{h}_0(\lambda^2,t)
\end{align*}

Then the second big term

\begin{align*}
    e^{i \lambda (1-x)}\Big[-\left[ G(\lambda, t) - G(-\lambda, t) \right]
+ i \lambda (e^{i \lambda}-e^{-i\lambda})\tilde{h}_0(\lambda^2,t)\Big]&= e^{i \lambda (1-x)}\Big[-\hat{u}_0(\lambda)+i\lambda\tilde{g}_0(\lambda^2,t)-e^{-i\lambda} i \lambda\tilde{h}(\lambda^2,t)\\
&+\hat{u}_0(-\lambda)+i\lambda \tilde{g}_0(\lambda^2,t)-e^{i\lambda}i\lambda \tilde{h}_0(\lambda^2,t)\\
&+i \lambda (e^{i \lambda}-e^{-i\lambda})\tilde{h}_0(\lambda^2,t)\Big]\\
&=e^{i \lambda (1-x)}\Big[-\hat{u}_0(\lambda)+\hat{u}_0(-\lambda)+2i\lambda\tilde{g}_0(\lambda^2,t)\\
&-2i\lambda e^{-i\lambda} \tilde{h}_0(\lambda^2,t)\Big]
\end{align*}

Putting it all together:
\begin{equation}\label{17}
\begin{aligned}
    u(x,t)&=\frac{1}{2\pi} \int_{-\infty}^{\infty} e^{i \lambda x
    -\lambda^2t} \hat{u}_0(\lambda)\:d\lambda\\
    &-\frac{1}{2\pi}\int_{\partial D^{+}} \frac{e^{ -\lambda^2t}}{\Delta \lambda}  \Bigg\{e^{i\lambda(x+1)}\hat{u}_0(\lambda)-e^{i\lambda(x-1)} \hat{u}_0(-\lambda)-2i\lambda e^{i \lambda(x-1)}\tilde{g}_0(\lambda^2,t)+2i\lambda e^{i\lambda x}\tilde{h}_0(\lambda^2,t)\\ 
&-e^{i \lambda (1-x)}\hat{u}_0(\lambda)+e^{i \lambda (1-x)}\hat{u}_0(-\lambda)+2i\lambda e^{i \lambda (1-x)}\tilde{g}_0(\lambda^2,t)-2i\lambda e^{-i\lambda x} \tilde{h}_0(\lambda^2,t)\Bigg\} \:d\lambda 
\end{aligned}
\end{equation}

Now, we will use the identities:
\begin{itemize}
    \item \(e^{i\lambda(x+1)}-e^{i\lambda(1-x)}=2i e^{i\lambda}\sin{(\lambda x)}\)
    \item \(-e^{i\lambda(x-1)}+e^{i\lambda(1-x)}=-e^{-i\lambda(1-x)}+e^{i\lambda(1-x)}=2i\sin{(\lambda(1-x))}\)
\end{itemize}

Simplifying the \(\eqref{17}\) we find

\begin{equation}\label{18}
\begin{aligned}
    u(x,t)&=\frac{1}{2\pi} \int_{-\infty}^{\infty} e^{i \lambda x
    -\lambda^2t} \hat{u}_0(\lambda)\,d\lambda\\
    &-\frac{1}{2\pi}\int_{\partial D^{+}} \frac{e^{ -\lambda^2t}}{\Delta (\lambda)}  \Bigg\{2\sin{(\lambda x)\Big[i e^{i\lambda}\hat{u}_0(\lambda)-2\lambda\tilde{h}_0(\lambda^2,t)\Big]}\\
    &+2\sin{(\lambda(1-x))}\Big[i\hat{u}_0(-\lambda)-2 \lambda\tilde{g}_0(\lambda^2,t)\Big]\Bigg\} \,d\lambda 
\end{aligned}
 \end{equation}

\section{Numerical implementation}\label{3}
We now specialise the general representation derived in Section~\ref{2} to the explicit choice of data
\begin{align*}
    u_0(x)=e^{-x},\quad g_0(t)=\cos{(t)},\quad h_0(t)=e^{-1} \cos{(t)}
\end{align*}

One verifies immediately that the compatibility conditions are satisfied:

\begin{align*}
    u(0,0)=1=u_0(0) \quad \text{and} \quad u(1,0)=e^{-1}=u_0(1) 
\end{align*}

thus the initial and boundary conditions are compatible both at \(x=0,t=0\) and \(x=1,t=0\). The transforms that correspond to these conditions are:

    \begin{align*}
\hat{u}_0(\lambda)
&=\int_0^1 e^{- x}e^{-i\lambda x}\,dx
=\int_0^1 e^{-(1 +i\lambda)x}\,dx
=\frac{1-e^{-(1 +i\lambda)}}{1 +i\lambda}=\frac{-i (1-e^{-(1+i\lambda)})}{\lambda-i}
\end{align*}

\begin{align*}
\tilde{g}_0(\lambda,t) &= \frac{1}{2} \int_0^t \left(e^{(\lambda + i)s} + e^{(\lambda - i)s} \right) ds = \frac{1}{2} \left( \frac{e^{(\lambda+i)t} - 1}{\lambda + i} + \frac{e^{(\lambda - i)t} - 1}{\lambda - i} \right)\\
&=\frac{1}{2}\Bigg(\frac{e^{(\lambda+i)t}}{\lambda+i}+ \frac{e^{(\lambda-i)t}}{\lambda-i}-\frac{1}{\lambda+i}-\frac{1}{\lambda-i}\Bigg)\\
&=\frac{1}{2}\Bigg(\frac{e^{\lambda t}\Big(e^{it}(\lambda-i)+e^{-it}(\lambda+i)\Big)}{\lambda^2+1}-\frac{2\lambda}{\lambda^2+1}\Bigg)\\
&=\frac{e^{\lambda t}(\lambda \cos{t}+\sin{t})-\lambda}{\lambda^2+1}
\end{align*}

and
\begin{align*}
    \tilde{h}_0(t)=\frac{\tilde{g}_0(\lambda,t)}{e}
\end{align*}

The first integral of the  \(\eqref{18}\) can now be written  as:
\begin{align*}
    &-\frac{i}{2\pi} \int_{\partial D^{+}} e^{i \lambda x-\lambda^2 t} \frac{1}{\lambda-i}\, d\lambda -\frac{i}{2\pi} \int_{\partial D^{-}} e^{i \lambda x-\lambda^2 t} \frac{e^{-1}e^{-i\lambda }}{\lambda-i}\, d\lambda\\
    &-\frac{i}{2\pi} \int_{\partial D^{+}} e^{-\lambda^2 t} \Big( \frac{e^{i\lambda x}}{\lambda-i}+  \frac{ e^{-1}e^{i\lambda(1-x)}}{\lambda+i}\Big) \, d\lambda
\end{align*}

We deform \(\partial D^+\) to the contour \(C^+\) to avoid the singularities
$\lambda_1=e^{i\frac{\pi}{4}}$ and $\lambda_2=e^{3i\frac{\pi}{4}}$ .

The contour \(C^+\)
 is taken to be the trapezoidal path defined by

\begin{equation}\label{eq_para}
    z = 
    \begin{cases} 
       z_1 (r)= (-l+\frac{3}{4}i)+re^{i\frac{5\pi}{6}} ,  & r \in (-\infty, -\ell), \\
       z_2 (s)= s + hi ,                            & s \in [-\ell,\, \ell],      \\
              z_3 (r)= (l+\frac{3}{4}i)+re^{i\frac{\pi}{6}} , &  r \in (\ell, +\infty)
    \end{cases}
\end{equation}
 with $\ell=\frac{3}{4}$. Substituting these expressions into the closed-form representation  \(\eqref{18}\) and introducing the auxiliary quantity

\[Q(\lambda,x)=-4\lambda \Bigg\{ \frac{\sin{(\lambda x)}}{e}+\sin{(\lambda(1-x))}\bigg\}\]

the solution may be written as

\begin{equation}
    \begin{aligned}
        &u(x,t)= \frac{1}{2 \pi} \int_{C^{+}} \frac{e^{-\lambda^2 t}}{e^{i \lambda}-e^{-i\lambda}} \frac{Q(\lambda,x)\cdot(\lambda^2-1)}{(\lambda^2+1)\cdot(\lambda^4+1)}\ d\lambda\\
        &-\frac{1}{2\pi} \int_{C^+} \frac{1}{e^{i \lambda}-e^{-i\lambda}} \frac{Q(\lambda,x)}{(\lambda^4+1)}\cdot [\lambda^2 \cos{(t)}+\sin{(t)}]\ d\lambda
    \end{aligned}
\end{equation}

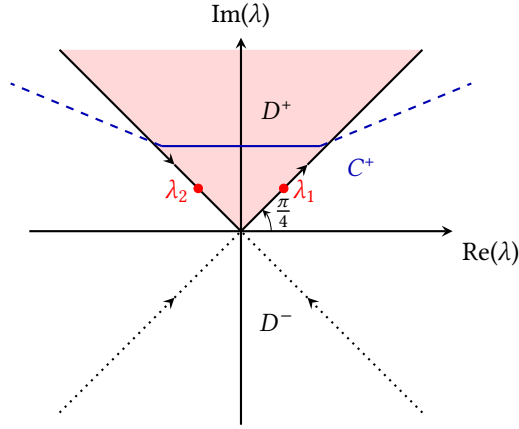
\begin{figure}[t]
\centering
\begin{tikzpicture}[scale=0.8, >=stealth]

\fill[red!15] (-3,3) -- (0,0) -- (3,3) -- cycle;

\draw[->, thick] (-3.5, 0) -- (3.5, 0) node[below right] {\(\text{Re}(\lambda)\)};
\draw[->, thick] (0, -3.2) -- (0, 3.2) node[above] {\(\text{Im}(\lambda)\)};

\draw[thick, ->] (-3,3) -- (-1.1,1.1);
\draw[thick] (-1.1,1.1) -- (0,0);
\draw[thick, ->] (0,0) -- (1.1,1.1);
\draw[thick] (1.1,1.1) -- (3,3);

\draw[dotted, thick, ->] (-3,-3) -- (-1.1,-1.1);
\draw[dotted, thick] (-1.1,-1.1) -- (0,0);
\draw[dotted, thick, ->] (3,-3) -- (1.1,-1.1);
\draw[dotted, thick] (1.1,-1.1) -- (0,0);

\node at (0.6,2) {\(D^+\)};
\node at (0.6,-1.5) {\(D^-\)};


\coordinate (l1) at ({1/sqrt(2)}, {1/sqrt(2)});
\coordinate (l2) at ({-1/sqrt(2)}, {1/sqrt(2)});
\filldraw[red] (l1) circle (2pt);
\filldraw[red] (l2) circle (2pt);
\node[anchor=north west, red] at ({1/sqrt(2)}, {0.3+1/sqrt(2)}) {\(\lambda_1\)};
\node[anchor=north east, red] at ({-1/sqrt(2)}, {0.3+1/sqrt(2)}) {\(\lambda_2\)};

\coordinate (bLeft) at ({-1.85/sqrt(2)}, {1/sqrt(2) + 0.7});
\coordinate (bRight) at ({1.85/sqrt(2)}, {1/sqrt(2) + 0.7});
\draw[thick, blue!70!black] (bLeft) -- (bRight);

\draw[thick, dashed, blue!70!black] (bLeft) -- ++(-2.5, {2.5 * 0.4142});
\draw[thick, dashed, blue!70!black] (bRight) -- ++(2.5, {2.5 * 0.4142});

\node[blue!70!black] at (2, {1.5/sqrt(2)}) {\(C^+\)};

\draw[->] (0.5,0) arc[start angle=0, end angle=45, radius=0.5];
  \node at (0.7, 0.3) {$\frac{\pi}{4}$};

\end{tikzpicture}
\caption{The trapezoidal contour $C^+$ used for numerical quadrature in \textsc{Maple}. The path is constructed to bypass the poles (red dots) while the linear segments are oriented to ensure the maximum exponential decay of the integrand $e^{-\lambda^2 t}$. } \label{fig_ex1b}
\end{figure}

\begin{figure}[t]
    \centering
    \includegraphics[scale=0.6]{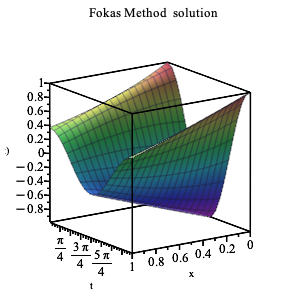}
   \caption{The solution \(u(x,t)\) in the range \(x \in [0,1]\) and \(t \in [0,2\pi]\) is presented a three-dimensional perspective plot.
}
\label{fig:3dex2}
    \end{figure}

We consider the integral component \[u_2= \frac{1}{2\pi} \int_{C^+} \frac{1}{e^{i \lambda}-e^{-i\lambda}} \frac{Q(\lambda,x)}{(\lambda^4+1)}\cdot [\lambda^2 \cos{(t)}+\sin{(t)}]\ d\lambda] d\lambda\]

This way, the roots are outside $C^+$, the integrand is analytic, and using Cauchy's theorem, we get
\begin{equation*}\label{eq_pole}
u_2(x,t)=0.    
\end{equation*}

The solution reduces to a single integral:
\[
u(x,t)=\frac{1}{2 \pi} \int_{C^{+}} \frac{e^{-\lambda^2 t}}{e^{i \lambda}-e^{-i\lambda}} \frac{Q(\lambda,x)\cdot(\lambda^2-1)}{(\lambda^2+1)\cdot(\lambda^4+1)}\ d\lambda
\]

  A three-dimensional visualization of the analytical solution \(u(x,t)\) for \(x \in[0,1]\)  and $t\in [0,2\pi]$  is shown in Figure~\ref{fig:3dex2}.

The numerical implementation of the Fokas Method using a fully deformed trapezoidal contour was carried out in Maple on a standard desktop machine. The computation required approximately 193.44~MiB of memory, with a CPU time of 11.25 seconds and a real (wall clock) time of 10.92 seconds.  

\section{Conclusion}
In this work, we applied the Unified Transform Method (Fokas method) to the one-dimensional heat equation on the finite interval [0,1]
, subject to Dirichlet boundary conditions. Starting from the spectral formulation, we derived the Global Relation and used it to eliminate the unknown boundary data \(\tilde{g}_1\)
and \(\tilde{h}_1\)
, arriving at a closed-form integral representation of the solution entirely in terms of prescribed data.
For an explicit example with exponential initial data \(u_0(x) = e^{-x}\)
and cosine boundary conditions \(g_0(t) = \cos(t)\), \(h_0(t) = e^{-1}\cos(t)\), the integral representation was simplified analytically. A key observation was that the contribution \(u_2\)
 vanishes by Cauchy's theorem, since the contour \(C^+\)
 is constructed to lie above the poles of the integrand, leaving a single, well-defined contour integral that was evaluated numerically in \textsc{Maple} using a trapezoidal contour parametrization.
The numerical computation confirmed the efficacy of the method, producing a smooth three-dimensional solution surface over \(x \in [0,1]\)
 and \(t \in [0, 2\pi]\)
 consistent with the prescribed initial and boundary data. The implementation required modest computational resources, completing in under 12 seconds.
This paper aims to serve as an accessible entry point to the Fokas method for students with a background in Fourier analysis and complex variables, demonstrating how contour deformation and the Global Relation combine to yield a practical solution scheme. Future directions include extending this approach to more complex geometries, systems of evolution equations, and applications to diffusion modeling in precision agriculture, where accurate heat and moisture transport models are critical for irrigation optimization.


\bibliographystyle{plain}
\bibliography{sample-base}

\appendix

\end{document}